\documentclass[12pt]{article}
\usepackage{mathrsfs}
\usepackage{}
\usepackage{amssymb}
\usepackage{amsfonts}
\usepackage{amsfonts}
\usepackage{latexsym}
\usepackage{CJK}
\date{}

\usepackage{amsmath}
\usepackage{amscd}
\usepackage{amsthm}
\usepackage{latexsym}
\usepackage[all]{xy}

\numberwithin{equation}{section}

\newcommand {\Dim}{\textrm{dim}}
\newcommand {\End}{\textrm{End}}
\newcommand {\Hom}{\textrm{Hom}}
\newcommand {\Ext}{\textrm{Ext}}
\newcommand {\Mod}{\textrm{mod}}

\newcommand{\bsm}{\begin{smallmatrix}}
\newcommand{\esm}{\end{smallmatrix}}

\begin{document}

\linespread{1.2}

\title{\bf Mutation graphs of maximal rigid modules over finite dimensional preprojective algebras$^\star$}
\author{{\small Hongbo Yin, Shunhua Zhang$^*$}\\
{\small  School of Mathematics,\ Shandong University,\ Jinan 250100,
P. R. China }}

\date{}
\maketitle

\pagenumbering{arabic}

\begin{center}
 \begin{minipage}{120mm}
   \small\rm
   {\bf  Abstract}\ \  Let $Q$ be a finite quiver of Dynkin type and
$\Lambda=\Lambda_Q$ be the  preprojective algebra of $Q$ over an
algebraically closed field $k$. Let $\mathcal {T}_\Lambda$ be the
mutation graph of maximal rigid $\Lambda$ modules.  Geiss, Leclerc
and Schr$\ddot{\rm o}$er conjectured that $\mathcal {T}_\Lambda$ is
connected, see [C.Geiss, B.Leclerc, J.Schr\"{o}er, Rigid modules
over preprojective algebras, Invent.Math., 165(2006), 589-632].   In
this paper, we prove that this conjecture is true when  $\Lambda$ is
of representation finite type or tame type. Moreover, we also prove
that $\mathcal {T}_\Lambda$ is isomorphic to the tilting graph of
${\rm End}_\Lambda T$ for each  maximal rigid $\Lambda$-module  $T$
if $\Lambda$ is representation-finite.
\end{minipage}
\end{center}

\vskip0.1in


\vskip0.2in


\footnote{ $^\star$Supported by the NSF of China (Grant No.
11171183).}

\footnote{ $^*$Corresponding author.}

\footnote{  Email addresses: \ yinhongbo0218@126.com(H.Yin),
 \ shzhang@sdu.edu.cn(S.Zhang).}

\vskip0.2in

\section{Introduction}

Let $Q$ be a finite quiver without oriented cycles and $kQ$ be the
path algebra  of $Q$ over an algebraically closed field $k$. The
preprojective algebra $\Lambda=\Lambda_Q$ of $Q$ was introduced by
Gelfand and Ponomarev in \cite{GP} such that $\Lambda$ contains $kQ$
as a subalgebra, and when considered as a left $kQ$ module,
$\Lambda$  decomposes as a direct sum of the indecomposable
preprojective $kQ$ modules with one from  each isomorphism class.
Now,  preprojective algebras play important roles in representation
theory and other areas of mathematics, such as resolutions of
Kleinian singularities, quantum groups, quiver varieties, and
cluster theory, see $\cite{CB, GLS1, GLS2, GLS3, GLS4, K, L}$ for
details.

\vskip0.2in

By using mutations of maximal rigid modules and their endomorphism
algebras over preprojective algebras of Dynkin type, Geiss, Leclerc
and Schr$\ddot{\rm o}$er studied the cluster algebra structure on
the ring $\mathbb{C}[N]$ of polynomial functions on a maximal
unipotent subgroup $N$ of a complex Lie group of Dynkin type, and
obtained that all cluster monomials of $\mathbb{C}[N]$ belong to the
dual semicanonical basis, see \cite{GLS1}.

\vskip0.2in

Let $Q$ be a Dynkin quiver, and $\Lambda$ be the preprojective
algebra of $Q$. Recall from \cite{GLS1}, $\mathcal {T}_\Lambda$
denotes  the \emph{mutation graph} of maximal rigid modules of
$\Lambda$. Fix a basic maximal rigid $\Lambda$-module $T$, then the
contravariant functor $F^T = \Hom_\Lambda(-, T) :  {\rm mod}\
\Lambda\rightarrow {\rm mod}\ {\rm End}_\Lambda  T$ yields an
anti-equivalence of categories $${\rm mod}\ \Lambda\rightarrow
\mathcal {P}({\rm mod}\ {\rm End}_\Lambda T)$$ where $\mathcal
{P}({\rm mod}\ {\rm End}_\Lambda T)\subset {\rm mod}\ {\rm
End}_\Lambda  T $  denotes the full subcategory of all ${\rm
End}_\Lambda  T$-modules of projective dimension at most one.
Moreover, the functor $F^T$ induces an embedding of graphs $\psi_T:
\mathcal {T}_\Lambda \rightarrow \mathcal {T}_{{\rm End}_\Lambda\
T}$ whose image is a union of connected components of $\mathcal
{T}_{{\rm End}_\Lambda\ T}$, where $\mathcal{T}_{{\rm End}_\Lambda\
T}$ is the tilting graph of the algebra ${\rm End}_\Lambda\ T$. Each
vertex of $\mathcal {T}_\Lambda $ (and therefore each vertex of the
image of $\psi_T$) has exactly $r-n$ neighbours.

\vskip0.2in

In \cite{GLS1}, Geiss, Leclerc and Schr$\ddot{\rm o}$er conjectured
that the graph $\mathcal {T}_\Lambda$ is connected. In this paper,
we prove that this conjecture is true when $\Lambda$ is of
representation finite type or tame type. Moreover, we also prove
that $\psi_T$ is an isomorphism whenever $\Lambda$ is representation
finite. The following theorems are our main results.

\vskip0.2in

{\bf Theorem 1.}\ {\it Let $\Lambda$ be a
preprojective algebra of type $A_n$ with $n\leq 4$, and  $T$ be a
maximal rigid $\Lambda$-module. Then the functor $F^T =
\Hom_\Lambda(-, T)$ induces an isomorphism of graphs $\psi_T:
\mathcal {T}_\Lambda \rightarrow \mathcal {T}_{{\rm End}_\Lambda\
T}$.}

\vskip0.2in

{\bf Corollary 2.}\ {\it Let $\Lambda$ be a
preprojective algebra of type $A_n$ with $n\leq 4$. Then for each
 maximal rigid $\Lambda$-module $T$, the tilting graph $\mathcal {T}_{{\rm End}_\Lambda\
T}$ of $\End_{\Lambda}T$ is isomorphic to the mutation graph
$\mathcal {T}_\Lambda$ of maximal rigid modules of $\Lambda$.}

\vskip0.2in

{\bf Remarks.}\  \  Let $\Lambda$ be a preprojective algebra with
finite representation type. The above  corollary implies that for
all maximal rigid $\Lambda$-modules, their endomorphism algebras
have same tilting graphs up to isomorphism.  However,  this kind of
algebras are very different, such as some of them is strongly
quasi-hereditary and most of them is even not quasi-hereditary, see
\cite{GLS4} for details.

\vskip0.2in

{\bf Theorem 3.}\ {\it Let $\Lambda$ be a preprojective algebra of
representation finite or tame type. Then the mutation graph
$\mathcal {T}_\Lambda$ of the maximal rigid $\Lambda$-modules is
connected.}

\vskip0.2in

This paper is organized as follows: in Section 2, we recall some
definitions and facts needed for our research, in Section 3, we
prove Theorem 1 and Corollary 2, in Section 4, we prove Theorem 3.

\vskip0.2in

\section{Preliminaries}

Let $k$ be  an algebraically closed field, and let A be a finite
dimensional algebra over $k$. We denote by ${\rm mod}~A$ the
category of all finitely generated left $A$-modules, and by ${\rm
ind}~A$ the full subcategory of ${\rm mod}~A$ consisting of one
representative from each isomorphism class of indecomposable
modules. For a $A$-module $M$, we denote by ${\rm add}~ M$ the full
subcategory of ${\rm mod}~A$ whose objects are the direct summands
of finite direct sums of copies of $M$. The projective dimension of
$M$ is denoted by ${\rm pd}~M$, and the Auslander Reiten translation
of $A$ by $\tau_A$.

\vskip0.2in

$T\in \Mod\ A$ is called a classical tilting module if
the following conditions are satisfied:\\
{(1)} ${\rm pd}~ T \leq 1$;\\
{(2)} ${\rm Ext}_A^{1}(T,T)=0 $;\\
{(3)} There is  an exact sequence
 $0\longrightarrow A \longrightarrow T_{0}\longrightarrow T_{1}
\longrightarrow 0$ with $T_{i}\in {\rm add}\ T $ for $0\leq i\leq
1$.

\vskip0.2in

Let ${\cal T}_{A}$ be the set of all basic classical tilting
$A$-modules up to isomorphism. According to \cite{GLS1, HU}, the
tilting graph $\mathcal {T}_A$  is the defined as following: the
vertices are the non-isomorphic basic tilting moduels, there is an
edge between $T_1$ and $T_2$ if $T_1=T'\oplus T_1'$ and
$T_2=T'\oplus T_2'$ for some $A$-module $T'$ and some indecomposable
$A$-modules $T'_1$ and $T'_2$ with $T'_1\not\simeq T'_2$.

\vskip0.2in

Let $Q=(Q_0, Q_1)$ be a connected quiver, where $Q_0$ is the set of
vertices and $Q_1$ is the set of arrows. Given an arrow $\alpha$, we
denote by $s(\alpha)$ the starting vertex of $\alpha$ and by
$t(\alpha)$ the ending vertex of $\alpha$. Let $\overline{Q}$ be the
double quiver of $Q$, which is obtained from $Q$ by adding an arrow
$\alpha^\ast: j\rightarrow i$ whenever there is an arrow $\alpha:
i\rightarrow j$ in $Q_1$. Let $Q_1^*=\{\alpha^*| \alpha\in Q_1\}$
and $\overline{Q}_1=Q_1\cup Q_1^*$. The preprojective algebra of $Q$
is defined as
\[\Lambda=\Lambda_Q=k\overline{Q}/(\rho)\]
where $\rho$ is the relation with
$$\rho=\sum_{\alpha\in Q_1}[\alpha, \alpha^*],$$
and $k\overline{Q}$ is the path algebra of $\overline{Q}$. See \cite{R1}.

\vskip0.2in

Note that the preprojective algebra $\Lambda$ is independent of the
orientation of $Q$, and that $\Lambda$ is finite dimensional if and
only if $Q$ is a Dynkin quiver.  Moreover, $\Lambda$ is also
self-injective if it is finite dimensional.  In particular,
$\Lambda$ is of finite representation type if and only if $Q$ is of
type $A_n$ with $n\leq 4$, and it is of tame representation type if
and only if $Q$ is of type $A_5$ or $D_4$, see \cite{DR, GLS2}.

\vskip0.2in

Let $d,e\in \mathbb{Z}^n$ be two dimension vectors.  The symmetry
bilinear form  is defined as $(d,e)=2\sum_{i\in
Q_0}d_ie_i-\sum_{a\in \overline{Q}_1}d_{s(a)}e_{t(a)}$. The
following lemma is proved in \cite{CB}.

\vskip0.2in

{\bf Lemma 2.1.}\ {\it Let $\Lambda$ be a  preprojective algebra and
$X, Y$ be $\Lambda$-modules. Then we have
 $${\rm dim~Ext}_{\Lambda}^1(X,Y)={\rm dim~Hom}_{\Lambda}(X,Y)+
 {\rm dim~Hom}_{\Lambda}(Y,X)-(\underline{\rm dim}~X,
 \underline{\rm dim}~Y).$$

In particular,    ${\rm dim~Ext}_{\Lambda}^1(X,Y)={\rm
dim~Ext}_{\Lambda}^1(Y,X)$.}

\vskip0.2in

From now on, we always assume that  $\Lambda$ is  a preprojective
algebra of Dynkin type. A $\Lambda$-module $T$ is called {\bf rigid}
if $\Ext_{\Lambda}^1(T,T)=0$. $T$ is called {\bf Maximal rigid} if
for any $\Lambda$-module $M$ with $\Ext_{\Lambda}^1(T\oplus
M,T\oplus M)=0$, then we have $M\in {\rm add}\ T$.

\vskip0.2in

Note that each maximal rigid $\Lambda$-module $T$ is also a
generator-cogenerator. Let $F^T=\Hom_{\Lambda}(-,T)$. A short exact
sequence $0\rightarrow X\rightarrow E\rightarrow Y\rightarrow 0$ of
$\Lambda$-modules is called $F^T$-exact if $ 0\rightarrow
F^T(Y)\rightarrow F^T(E)\rightarrow F^T(X)\rightarrow 0$ is an exact
sequence of ${\rm End}_\Lambda T $-modules. We denote by $F^T(Y,X)$
the equivalent classes of all the $F^T$-exact sequences as above.

\vskip0.2in

Let $\chi_T$ be a subcategory of $\Mod\ \Lambda$ whose objects admit
an ${\rm add}\ T$-resolution. Namely,  $X\in \chi_T$ if and only if
there is an exact sequence
$$\xymatrix{0\ar[r]&X\ar[r]&T_0\ar[r]&T_1\ar[r]&T_2\ar[r]&\cdots}$$
with all $T_i\in {\rm add}~T$, which  is  still exact by applying
the functor $\Hom_{\Lambda}(T,-)$. Let $\Ext_{F^T}^i(Y,X)$ be the
cohomology group by applying the functor $\Hom_{\Lambda}(Y,-)$ to an
${\rm add}\ T$-resolution of $X$.

\vskip0.2in

The following lemma is proved in \cite{AS1, AS2}.

\vskip0.2in

{\bf Lemma 2.2.}\ {\it  Assume that $X\in \chi_T$ and $Y\in \Mod\
\Lambda$. Then  there are following functorial isomorphisms:

\vskip0.1in

(1)\ ${\rm Ext}_{F^T}^1(Y,X)\cong F^T(Y,X);$

\vskip0.1in

(2)\ ${\rm Ext}_{F^T}^i(Y,X)\cong {\rm Ext}_{\rm
End_{\Lambda}(T)}^i({\rm Hom}_{\Lambda}(X,T),{\rm
Hom}_{\Lambda}(Y,T))$ for all $i\geq 1$}.

\vskip0.2in

Let $\Lambda$ be a finite dimensional preprojective algebra,  and
let $T$ be a maximal rigid $\Lambda$-module. Then $\chi_T=\Mod\
\Lambda$ since every $\Lambda$-module has an add\ $T$-resolution
\cite[Corollary 5.2]{GLS1}.

\vskip0.2in

Recall from \cite[section 6]{GLS1}, the \emph{mutation graph}
$\mathcal {T}_\Lambda$ of maximal rigid modules is defined as
following.  The  vertex set of $\mathcal {T}_\Lambda$ is the set of
the isomorphism classes of basic maximal rigid $\Lambda$-modules,
and there is an edge between vertices $T_1$ and $T_2$ if and only if
$T_1 = T \oplus T'_1$ and $T_2 = T\oplus T'_2$ for some $T$ and some
indecomposable modules $T'_1$ and $T'_2$ with $T'_1\not\simeq T'_2$.

\vskip0.2in

{\bf Lemma 2.3.}\ {\it  Let $T$ be a basic maximal rigid
$\Lambda$-module. The functor $F^T: {\rm mod}~\Lambda \rightarrow
{\rm mod}~{\rm End}_{\Lambda}(T)$ induces an embedding of graphs
$\psi_T:\mathcal {T}_\Lambda\rightarrow \mathcal {T}_{{\rm
End}_{\Lambda}(T)}$ whose image is a union of connected components
of $\mathcal {T}_{{\rm End}_{\Lambda}(T)}$.}

\vskip0.2in

We follow the standard terminology and notation used in the
representation theory of algebras, see \cite{ASS, ARS, R}.

\section{The mutation graph and the tilting graph of representation finite preprojective algebras}

In this section, we assume that $\Lambda$ is a preprojective algebra
of representation finite type. Namely, $\Lambda$ is of type $A_n$
with  $n\leq 4$. For the AR-quivers of this kind of preprojective
algebras we refer to  \cite[section 20.1]{GLS2}. Here we give the
stable AR-quivers of $\Lambda_{A_3}$ and $\Lambda_{A_4}$ for
convenience.

\vskip0.2in

\[\xymatrix@C=0.3cm@R=0.3cm{
          \ar@{.}[d]&\circ\ar[dr]&&\circ\ar[dr]&&Y\ar[dr]&\ar@{.}[d]\\
          \circ\ar[ur]\ar[rd]\ar@{.}[d]&&Z_1\ar[ur]\ar[dr]&&\circ\ar[ur]\ar[dr]&&Z_2\ar@{.}[d]\\
          &X\ar[ur]&&\circ\ar[ur]&&Z_3\ar[ur]&\\
           }\]
\begin{center} The stable quiver of $\Lambda_{A_3}$  \end{center}

\[\xymatrix@C=0.3cm@R=0.3cm{\circ\ar@{.}[dd]\ar[rd]&&\circ\ar[rd]&
&\bullet\ar[rd]&&\circ\ar[rd]&&\circ\ar[rd]&&\circ\ar[rd]&&\circ\ar@{.}[dd]\\
       &\circ\ar[ru]\ar[rd]&&\bullet\ar[ru]\ar[rd]&&\circ\ar[ru]
       \ar[rd]&&\circ\ar[ru]\ar[rd]&&\circ\ar[ru]\ar[rd]&&\circ\ar[ru]\ar[rd]\\
         \circ\ar@{.}[d]\ar[ru]\ar[rd]&&\circ\ar[ru]\ar[rd]&
         &\bullet\ar[ru]\ar[rd]&&\circ\ar[ru]\ar[rd]&&\circ\ar[ru]\ar[rd]&&
         \circ\ar[ru]\ar[rd]&&\circ\ar@{.}[d]\\
         \circ\ar@{.}[d]\ar[r]&\circ\ar[ru]\ar[rd]\ar[r]&\circ\ar[r]&
         \bullet\ar[ru]\ar[rd]\ar[r]&\bullet\ar[r]&\circ\ar[ru]\ar[rd]\ar[r]&\circ\ar[r]&
         \circ\ar[ru]\ar[rd]\ar[r]&\circ\ar[r]&
         \circ\ar[ru]\ar[rd]\ar[r]&\circ\ar[r]&\circ\ar[ru]\ar[rd]\ar[r]&\circ\ar@{.}[d]\\
         \circ\ar[ru]&&\circ\ar[ru]&&\bullet\ar[ru]&&\circ\ar[ru]&&\circ\ar[ru]&&\circ\ar[ru]&&\circ}\]
\begin{center} The stable quiver of $\Lambda_{A_4}$  \end{center}

\vskip0.2in

{\bf Definition.}\ {\it Two AR-sequences are called {\bf centrally
connected} if they have common indecomposable summands in the middle
terms. A {\bf column} in the AR-quiver is a set consist of the
indecomposable summands of the middle terms in the centrally
connected AR-sequences.}

\vskip0.1in

 {\it A path from $X$ to $Y$ in the AR-quiver is a chain of irreducible morphisms
$X=M_0\rightarrow M_1\rightarrow M_2\rightarrow\cdots\rightarrow
M_{n-1}\rightarrow M_n=Y$. We say that $Z$ is between $X$ and $Y$ if
there is a chain
$$X=M_0\rightarrow M_1\rightarrow
M_2\rightarrow\cdots\rightarrow M_{n-1}\rightarrow M_n=Y$$ such that
all $M_i$ is not in the same column with $Y$ for $0<i<n$ and that
$Z$ is in the same column with some one of $M_i$ with $0<i<n$.}

\vskip0.1in

 {\it  A class $\Sigma$ of pairwise non-isomorphic indecomposable
 $\Lambda$-modules in the stable quiver above is
called a complete slice if it satisfies the following conditions:

\vskip0.1in

(1) the indecomposable modules in $\Sigma$ lie in different
$\tau$-orbits;

\vskip0.1in

(2) $\Sigma$ is convex. Namely, if $X$ and $Y$ belong to $\Sigma$
and there is a path from $X$ to $Z$ and a path from $Z$ to $Y$, then
$Z$ belongs to $\Sigma$.}

\vskip0.1in

 {\it  A complete slice is  called standard if it lies in two adjacent columns.}

\vskip0.2in

For example, in the stable quiver of $\Lambda_{A_3}$, $Z_1$ is
between $X$ and $Y$ while $Z_2$ is between $Y$ and $X$. The complete
slice which consists of $\bullet$ in the stable quiver of
$\Lambda_{A_4}$ is standard.

\vskip0.2in

{\bf Lemma 3.1.}\ {\it Given a communicative diagram of exact
sequences
$$\xymatrix{0\ar[r]&X\ar[r]^i\ar@{=}[d]&E\ar[r]\ar[d]&Y\ar[r]\ar[d]^h&0\\
          0\ar[r]&X\ar[r]^f&F\ar[r]&Z\ar[r]&0}$$
with the bottom sequence non-split. Then the top sequence is
non-split if and only if $h$ cannot factor through $F$.}

\vskip0.1in

{\bf Proof.}\  Apply the functor $\Hom_{\Lambda}(Y,-)$ to the bottom
 sequence, we get an exact sequence
$$
\xymatrix@C=0.5cm{0\ar[r]&\Hom_{\Lambda}(Y,X)\ar[r]&\Hom_{\Lambda}(Y,F)\ar[r]^{\alpha}&
\Hom_{\Lambda}(Y,Z)\ar[r]^{\beta}&\Ext_{\Lambda}^1(Y,X)}.
$$
Then $h$ is in the kernel of $\beta$ if and only if it is in the
image of $\alpha$. Namely, the top sequence is the zero element in
${\rm Ext}^1_{\Lambda}(Y, X)$ if and only if $h$ factors through
$F$. This complete the proof. $\hfill\Box$

\vskip0.2in

{\bf Lemma 3.2.}\ {\it  Let $\Lambda$ be a preprojective algebra of
type $A_n$, $n\leq4$. Let $X$, $Y$ and $Z$ be non-isomorphic
indecomposable $\Lambda$-modules with ${\rm Ext}_{\Lambda}^1(Y,
X)\neq 0$ and ${\rm Ext}_{\Lambda}^1(Z, X)\neq 0$. If $Y$ is between
$X$ and $Z$ with ${\rm Hom}_{\Lambda}(Y,Z)\neq 0$, then there is a
non-split exact sequence
$$(1)\ \ \ \ 0\rightarrow X\rightarrow E\rightarrow Y\rightarrow 0$$
which is induced from a non-split exact sequence
$$(2)\ \ \ \ 0\rightarrow X\stackrel{f}\longrightarrow F\rightarrow Z\rightarrow 0.$$}

\vskip0.1in

{\bf Proof.}\  Let (3) $0\rightarrow X\stackrel{i}\longrightarrow
M\rightarrow \tau^{-1}X\rightarrow 0$ be the AR-sequence start at
$X$. Then we have the following communicative diagram:
$$\xymatrix{0\ar[r]&X\ar[r]^i\ar@{=}[d]&M\ar[r]\ar[d]&\tau^{-1}X\ar[r]\ar[d]^h&0\\
          0\ar[r]&X\ar[r]^f&F\ar[r]&Z\ar[r]&0}.
$$
By using AR-formula ${\rm \overline{Hom}}_\Lambda(\tau^{-1}X,
Z)\simeq D{\rm Ext}^1_\Lambda(Z, X)$, we know that  different
sequences of the form (2) corresponds to different homomorphisms
from $\tau^{-1}X$ to $Z$ in the stable category
$\underline{\Mod}\Lambda$. According to Lemma 3.1, we know that $h$
can't factor through $F$.

Let $X$, $Y$ and $Z$ be non-isomorphic indecomposable
$\Lambda$-modules with $\Ext_{\Lambda}^1(Y, X)\neq 0\neq
\Ext_{\Lambda}^1(Z, X)$. If $Y$ is between $X$ and $Z$ with
$\Hom_{\Lambda}(Y,Z)\neq 0$, then by reading the pictures given in
\cite[section 20.4]{GLS2} we  know that there is a path from
$\tau^{-1}X$ to $Z$ which induces a nonzero morphism from
$\tau^{-1}X$ to $Z$ in $\underline{\Mod}\Lambda$ factoring through
$Y$. Hence there exists a morphism $g$ from $Y$ to $Z$ which cannot
factor through $F$. Then we have a pull-back diagram:
$$\xymatrix{0\ar[r]&X\ar[r] \ar@{=}[d]&E\ar[r]\ar[d]& Y\ar[r]\ar[d]^g&0\\
          0\ar[r]&X\ar[r]^f&F\ar[r]&Z\ar[r]&0},
$$
by Lemma 3.1 again, we have a non-split sequence of the form (1)
which is induced from (2).  $\hfill\Box$

\vskip0.2in

{\bf Remark.}  We should mention that Lemma 3.2 is not true without
the assumption that $Y$ is between $X$ and $Z$. The following
example is pointed out to us by C.M.Ringel. Let
\[\xymatrix{1\ar@<1ex>[r]^{\alpha_1}&2\ar@<1ex>[r]^{\alpha_2}
\ar@<1ex>[l]^{\alpha_1^*}&3\ar@<1ex>[r]^{\alpha_3}\ar@<1ex>[l]^{\alpha_2^*}&
            4\ar@<1ex>[l]^{\alpha_3^*}}\]
be the quiver of $\Lambda_{A_4}$. Take
$X=\begin{array}{c}4\\3\end{array}$, $Y=\begin{array}{c}2\\1\ \
3\\2\end{array}$, $Z=2$, $V=\begin{array}{c}2\ \ 4\\3\end{array}$.
Then $\Ext_{\Lambda}^1(Y,X)=\Ext_{\Lambda}^1(Z,X)=k$,
$\Hom_{\Lambda}(Y,Z)=k$. $0\rightarrow X\rightarrow P(2)\rightarrow
Y\rightarrow 0$ and $0\rightarrow X\rightarrow V\rightarrow
Z\rightarrow 0$ are the corresponding exact sequences. But the first
sequence cannot be induced by the second one because the inclusion
map $0\rightarrow X\rightarrow V$ cannot factor through $P(2)$,
since there is no map from $P(2)$ to the simple module $4$.

\vskip0.2in

{\bf Lemma 3.3.}\ {\it Let $\Lambda$ be a preprojective algebra of
type $A_n$, $n\leq4$. Let $X$ and $Y$ be non-isomorphic
indecomposable $\Lambda$-modules with ${\rm
Ext}_{\Lambda}^1(X,Y)\neq 0$. Let $N$ be an indecomposable
non-projective $\Lambda$-module which is between $X$ and $Y$ or in
the same column with $X$. Then any exact sequence
$$(*)\ \ \ \ \ \ \
0\rightarrow Y\rightarrow M\rightarrow X\rightarrow 0$$ is
$F^N$-exact. Moreover, if $N$ is in the same column with $X$, then
any exact sequence
$$(**)\ \ \ \ \ \ \
0\rightarrow X\rightarrow M\rightarrow Y\rightarrow 0$$ is also
$F^N$-exact. }

\vskip0.1in

{\bf Proof.}\  We choose a standard complete slice which contains
$X$ and extend it to a maximal rigid $\Lambda$ module $T$ by adding
all the indecomposable projective-injective modules. Then it follows
from the stable quiver of $\Lambda$ that every non-zero map from $Y$
to $N$ factors through $T$ since $Y\not\in\rm add~T$.

Note that $\Ext_{\Lambda}(X,T)=0$, by applying $\Hom_{\Lambda}(-,T)$
to the exact sequence ($\ast$), we get an exact sequence
$$0\rightarrow \Hom_{\Lambda}(X,T)\rightarrow
\Hom_{\Lambda}(M,T)\rightarrow \Hom_{\Lambda}(Y,T)\rightarrow 0.
$$
Thus every map from $Y$ to $T$ factors through $M$, which implies
that every map from $Y$ to $N$ factors through $M$. Namely, the
sequence
$$
0\rightarrow \Hom_{\Lambda}(X,N)\rightarrow
\Hom_{\Lambda}(M,N)\rightarrow \Hom_{\Lambda}(Y,N)\rightarrow 0
$$ is exact.

Now, we assume that $N$ is in the same column with $X$. Then any map
from $X$ to $N$ in the stable quiver factors through the maximal
rigid module obtained from the standard complete slice which
contains $X$. Repeat the proof above we see that any exact sequence
$$ (\ast\ast)\ \ \ \ \ \ \
0\rightarrow X\rightarrow M\rightarrow Y\rightarrow 0$$ is also
$F^N$-exact. This completes the proof. $\hfill\Box$

\vskip0.2in

{\bf Lemma 3.4.}\ {\it  Let $\Lambda$ be a preprojective algebra of
type $A_n$, $n\leq4$. Let $X$ and $Y$ be non-isomorphic
indecomposable $\Lambda$-modules with ${\rm
Ext}_{\Lambda}^1(X,Y)\neq 0$. Let  $N$ be an indecomposable
non-projective $\Lambda$-module. Then there exists a non-split exact
sequence
$$
(\ast)\ \ \ \ \ \ \ \ \ \  0\rightarrow Y\rightarrow E\rightarrow
X\rightarrow 0
$$ or
$$
(\ast\ast)\ \ \ \ \ \ \ \ \ \  0\rightarrow X\rightarrow
M\rightarrow Y\rightarrow 0
$$  which is $F^N$-exact. }

\vskip0.1in

{\bf Proof.}\ If $N$ is between $X$ and $Y$ or in the same column
with $X$, then any exact sequence ($\ast$) $0\rightarrow
Y\rightarrow M\rightarrow X\rightarrow 0$ is $F^N$-exact by Lemma
3.3. If $N$ is between $Y$ and $X$ or in the same column with $Y$,
then any exact sequence ($\ast\ast$) $0\rightarrow X\rightarrow
E\rightarrow Y\rightarrow 0$ is $F^N$-exact  by Lemma 3.3 again.
   $\hfill\Box$

\vskip0.2in

{\bf Lemma 3.5.}\ {\it  Let $\Lambda$ be a preprojective algebra of
type $A_n$, $n\leq4$. Let $X$ and $Y$ be non-isomorphic
indecomposable $\Lambda$-modules with ${\rm
Ext}_{\Lambda}^1(X,Y)\neq 0$. Let $N_1$ and $N_2$ be two
non-isomorphic indecomposable $\Lambda$-module with ${\rm
Ext}_{\Lambda}^1(N_1, N_2)= 0$. Then there exists a non-split exact
sequence
$$
(\ast)\ \ \ \ \ \ \ \ \ \  0\rightarrow Y\rightarrow E\rightarrow
X\rightarrow 0
$$ or
$$
(\ast\ast)\ \ \ \ \ \ \ \ \ \  0\rightarrow X\rightarrow
M\rightarrow Y\rightarrow 0
$$
which is both $F^{N_1}$-exact and $F^{N_2}$-exact.}

\vskip0.1in

{\bf Proof.}\ If $N_1$ and $N_2$ are both between $X$ and $Y$ or
both between $Y$ and $X$, then the assertion is true by Lemma 3.3.

If $N_1$ is in the same column with $X$, then both ($\ast$) and
($\ast\ast$) are exact by Lemma 3.3. Hence the assertion is  true by
Lemma 3.4.

Now, we assume that $N_1$ is between $X$ and $Y$ while $N_2$ is
between $Y$ and $X$.

If $\Hom_{\Lambda}(Y, N_2)=0$, then ($\ast$) is $F^{N_2}$-exact, and
by Lemma 3.3, ($\ast$) is also $F^{N_1}$-exact.

If $\Hom_{\Lambda}(Y, N_2)\neq 0$, then by Lemma 3.3, any exact
sequence of form ($\ast$) is $F^{N_1}$-exact and any exact sequence
of form ($\ast\ast$) is $F^{N_2}$-exact.

{\bf Case I.}\  If there exists a non-split sequence of form
($\ast\ast$) is $F^{N_1}$-exact, then our sequence is true.

{\bf Case II.}\   Now, we suppose that any  non-split exact sequence
of the form $(\ast\ast)$ is not $F^{N_1}$-exact. We claim that
$\Ext_{\Lambda}^1(N_2, X)=0$.

Indeed, if by contrary we assume that $\Ext_{\Lambda}^1(N_2, X)\neq
0$, then by Lemma 3.2, there exists a non-split exact sequence
$$0\rightarrow X\stackrel{j}\longrightarrow E\rightarrow
Y\rightarrow 0
$$ which is induced from a non-split  exact sequence
$0\rightarrow X\stackrel{f}\longrightarrow F\rightarrow
N_2\rightarrow 0$. Then we have the following commutative diagram:
$$\xymatrix{0\ar[r]&X\ar[r]^j\ar@{=}[d]&E\ar[r]\ar[d]&Y\ar[r]\ar[d]^h&0\\
          0\ar[r]&X\ar[r]^f&F\ar[r]&N_2\ar[r]&0} .$$
Thus $f$ factors through $j$.

Note that $0\rightarrow X\rightarrow E\rightarrow Y\rightarrow 0$ is
not $F^{N_1}$-exact, hence there exists a map $\lambda$ from $X$ to
$N_1$ which cannot factor through $E$, this forces that there exists
a map $g$ from $X$ to $N_1$ such that $g$ cannot factor through $F$.

Then we have  following push-out diagram:
$$\xymatrix{0\ar[r]&X\ar[r]^f\ar[d]^g&F\ar[r]\ar[d]&N_2\ar[r]\ar@{=}[d]&0\\
          0\ar[r]&N_1\ar[r]&M\ar[r]&N_2\ar[r]&0} ,$$
which implies that the exact sequence $0\rightarrow N_1\rightarrow
M\rightarrow N_2\rightarrow 0$ is non-split. This is a contradiction
with $\Ext^1_\Lambda(N_2, N_1)=\Ext^1_\Lambda(N_1, N_2)=0$. Hence
our claim is true. Namely, $\Ext^1_\Lambda(X,
N_2)=\Ext^1_\Lambda(N_2, X)=0$. Therefore $0\rightarrow Y\rightarrow
M\rightarrow X\rightarrow 0$ is $F^{N_1}$-exact and $F^{N_2}$-exact.
This completes the proof. $\hfill\Box$

\vskip0.2in

{\bf Lemma 3.6.}\ {\it Let $\Lambda$ be a preprojective algebra of
type $A_n$, $n\leq4$. Let $X$ and $Y$ be indecomposable
$\Lambda$-modules. Then there exits an dense open orbits in the
variety of extensions between $X$ and $Y$.}

{\bf Proof.}\ It can be proved easily from \cite[section
2.1]{B1}.$\hfill\Box$

\vskip0.2in

{\bf Remark.} Recall from \cite{B1}, we say that  $M$ degenerate to
$N$ and denote by $M\leq_{deg}N$, if
$\mathcal{O}_N\subset\overline{\mathcal{O}}_M$. Let $\Lambda$ be a
preprojective algebra of type $A_n$, $n\leq4$. Using the AR-formula
and hammock algorithm we can see that $\Dim~
\Ext_{\Lambda}^1(X,Y)\leq 2$. In the case of $A_3$, we have that
$\Dim ~\Ext_{\Lambda}^1(X,Y)\leq 1$. If $\Dim~
\Ext_{\Lambda}^1(X,Y)= 2$, by Lemma 3.6 we have two non-split exact
sequence
$$0\rightarrow Y\rightarrow M_1\rightarrow X\rightarrow 0$$
and
$$0\rightarrow Y\rightarrow M_2\rightarrow X\rightarrow 0$$
such that $M_1\leq_{deg} M_2$. Then by \cite{B1}, we know that
$$\Dim\ \Hom_{\Lambda} (M_1, T)\leq \Dim\ \Hom_{\Lambda}(M_2, T)$$ for
any $\Lambda$-module $T$. Hence if $$0\rightarrow Y\rightarrow
M_1\rightarrow X\rightarrow 0$$ is $F^T$-exact, then $$0\rightarrow
Y\rightarrow M_2\rightarrow X\rightarrow 0$$ is also $F^T$-exact by
comparing dimensions.

\vskip0.2in

{\bf Lemma 3.7.}\ {\it Let $\Lambda$ be a preprojective algebra of
type $A_n$ with $n\leq4$. Let $X$ and $Y$ be non-isomorphic
indecomposable $\Lambda$-modules with ${\rm
Ext}_{\Lambda}^1(X,Y)\neq 0$. Let $T$ be a basic maximal rigid
$\Lambda$-module. Then there exists a non-split exact sequence
$$
(\ast)\ \ \ \ \ \ \ \ \ \  0\rightarrow Y\rightarrow E\rightarrow
X\rightarrow 0
$$ or
$$
(\ast\ast)\ \ \ \ \ \ \ \ \ \  0\rightarrow X\rightarrow
M\rightarrow Y\rightarrow 0
$$
which is $F^{T}$-exact.}

{\bf Proof.}\ According to  the Remark after Lemma 3.6, we only need
to consider the case that $\Dim\ \Ext_{\Lambda}^1(X,Y)=1$.

If $\Lambda$ is of type $A_3$, then $T$ has three indecomposable non
projective direct summands $T_1$, $T_2$ and $T_3$. We divide them
into three combinations $\{T_1, T_2\}$, $\{T_2, T_3\}$ and $\{T_1,
T_3\}$. By Lemma 3.5, there is an exact sequence $(\ast)$ or
$(\ast\ast)$ which is $F^{T_i}$-exact for at least two combinations,
then it is $F^{T}$-exact.

If $\Lambda$ is of type $A_4$, then $T$ has six indecomposable non projective direct
summands $T_1$, $T_2$, $T_3$, $T_4$, $T_5$ and $T_6$. There are twenty combinations say
$\{C_i\}_{1\leq i\leq 20}$,
such that each $C_i$ consists of three non isomorphic direct summands. We say an exact
sequence is $C_i$-exact if it is $F^{T_k}$-exact with $T_k\in C_i$. Then as above each
$C_i$ has at least one exact sequence $(\ast)$ or $(\ast\ast)$ that is $C_i$-exact.

Now we show the assertion that if we cut the set $\{C_i\}$ into two
parts, there always exists one part that covers all the six $T_i$.
If we choose three elements from five elements, there is ten kind of
possibilities. So, if we cut $\{C_i\}$ into two parts $U_1$ and
$U_2$ such that the number of $C_i$ in $U_1$ is bigger than ten,
then $\bigcup_{C_i\in U_1}C_i$ must contain at least six elements.
The assertion is right. Now suppose the number of the $C_i$ in $U_1$
and $U_2$ are both ten. If $\bigcup_{C_i\in U_1}C_i$ contains six
elements, then the assertion is right. If not, $\bigcup_{C_i\in
U_1}C_i$ contains five elements. Without loss of generality, we may
assume that $\bigcup_{C_i\in U_1}C_i=\{T_1, T_2, T_3, T_4, T_5\}$.
Then each $C_i$ in $U_2$ contains $T_6$ and $\bigcup_{C_i\in
U_2}C_i=\{T_1, T_2, T_3, T_4, T_5, T_6\}$. The assertion is also
true.

Now we divide $\{C_i\}$ into two parts according to the $C_i$-exact
sequence is of form $(\ast)$ or of form $(\ast\ast)$. If $C_i$ can
belong to both part, put it in only one part. Then there is an exact
sequence that is $F^{T_i}$-exact for all $1\leq i\leq 6$.
$\hfill\Box$

\vskip0.2in {\bf Proposition 3.8.}\ {\it\  Let $\Lambda$ be a
preprojective algebra of type $A_n$, $n\leq4$. Let $T$ be a basic
maximal rigid $\Lambda$-module and $B={\rm End}~ T$. Then every
classical tilting $B$-module is of the form ${\rm Hom}_{\Lambda}(T',
T)$, where $T'$ is a maximal rigid $\Lambda$-module.}

\vskip0.1in

{\bf Proof.}\ By Proposition 4.4 in \cite{GLS1}, we know that any
$B$-module with projective dimension at most 1 is of the form
$\Hom_{\Lambda}(M, T)$ with $M$ being a $\Lambda$-module.

If $M$ is not rigid, then there are indecomposable direct summands
$X$ and $Y$ of $M$ such that
$\Ext_{\Lambda}^1(X,Y)=\Ext_{\Lambda}^1(Y,X)\neq 0$. By Lemma 2.2 and
Lemma 3.7, we know that
$$\Ext_{\End_{\Lambda}(T)}^1(\Hom_{\Lambda}(X,T),\Hom_{\Lambda}(Y,T))\neq
0$$ or
$$\Ext_{\End_{\Lambda}(T)}^1(\Hom_{\Lambda}(Y,T),\Hom_{\Lambda}(X,T))\neq
0.$$ Hence $\Hom_{\Lambda}(M, T)$ is not partial tilting as $B$-module.

In particular, any partial tilting $B$-module is of the form
$\Hom_{\Lambda}(T', T)$ with $T'$ being a rigid $\Lambda$-module.
Note that the number of non-isomorphic simple $B$-module is equal to
the number of the non-isomorphic indecomposable direct summands of
the maximal rigid $\Lambda$-module $T$, hence $\Hom_{\Lambda}(T',
T)$ is a tilting module if and only if $T'$ is a maximal rigid
$\Lambda$-module. This completes the proof. $\hfill\Box$

\vskip0.2in

Summarizing above discussions, we have the following theorem which
is one of our main results.

\vskip0.2in

{\bf Theorem 3.9.}\ {\it\ Let $\Lambda$ be a preprojective algebra
of type $A_n$ with $n\leq 4$, and  $T$ be a maximal rigid
$\Lambda$-module. Then the functor $F^T = {\rm Hom}_\Lambda(-, T)$
induces an isomorphism of graphs $\psi_T: \mathcal {T}_\Lambda
\rightarrow \mathcal {T}_{{\rm End}_\Lambda\ T}$.}

\vskip0.1in

{\bf Proof.}\ By Lemma 2.3, we know that $\psi_T$ is injective, and
$\psi_T$ is also surjective by Proposition 3.8. Namely, $\psi_T$ is
an isomorphism. This completes the proof. $\hfill\Box$

\vskip0.2in

The following Corollary is a direct consequence.

\vskip0.2in

{\bf Corollary 3.10.}\ {\it\ Let $\Lambda$ be a preprojective algebra
of type $A_n$ with $n\leq 4$. Then for each
 maximal rigid $\Lambda$-module $T$, the tilting graph $\mathcal {T}_{{\rm End}_\Lambda\
T}$ of $\End_{\Lambda}T$ is isomorphic to the mutation graph
$\mathcal {T}_\Lambda$ of maximal rigid modules of $\Lambda$.}

\vskip0.2in

We illustrate our results by the example  $\Lambda_{A_3}$. The
AR-quiver of $\Lambda_{A_3}$ is as follows, here we represent
$\Lambda_{A_3}$-modules by Lowvey series.

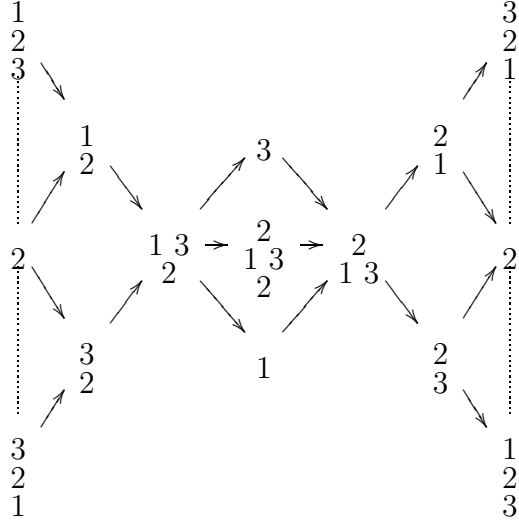
\renewcommand\arraycolsep{0.1cm}
\newcommand{\one}{{\begin{array}{c}\vspace{-10pt}1\\\vspace{-10pt}2\\\vspace{-10pt}3\end{array}}}
\newcommand{\two}{{\begin{array}{c}\vspace{-10pt}3\\\vspace{-10pt}2\\\vspace{-10pt}1\end{array}}}
\newcommand{\three}{{\begin{array}{c}\vspace{-10pt}1\\\vspace{-10pt}2\\\end{array}}}
\newcommand{\four}{{\begin{array}{c}\vspace{-10pt}3\\\end{array}}}
\newcommand{\five}{{\begin{array}{c}\vspace{-10pt}2\\\vspace{-10pt}1\\\end{array}}}
\newcommand{\six}{{\begin{array}{c}\vspace{-10pt}2\\\end{array}}}
\newcommand{\seven}{{\begin{array}{c}\vspace{-10pt}1\ 3\\\vspace{-10pt}2\\\end{array}}}
\newcommand{\eight}{{\begin{array}{c}\vspace{-10pt}2\\\vspace{-10pt}1\ 3\\\vspace{-10pt}2\end{array}}}
\newcommand{\nine}{{\begin{array}{c}\vspace{-10pt}2\\\vspace{-10pt}1\ 3\\\end{array}}}
\newcommand{\ten}{{\begin{array}{c}\vspace{-10pt}2\\\end{array}}}
\newcommand{\eleven}{{\begin{array}{c}\vspace{-10pt}3\\\vspace{-10pt}2\\\end{array}}}
\newcommand{\twive}{{\begin{array}{c}\vspace{-10pt}1\\ \end{array}}}
\newcommand{\thirteen}{{\begin{array}{c}\vspace{-10pt}2\\\vspace{-10pt}3\\\end{array}}}
\newcommand{\fourteen}{{\begin{array}{c}\vspace{-10pt}3\\\vspace{-10pt}2\\\vspace{-10pt}1\end{array}}}
\newcommand{\fifteen}{{\begin{array}{c}\vspace{-10pt}1\\\vspace{-10pt}2\\\vspace{-10pt}3\end{array}}}

\[\xymatrix@C=0.3cm@R=0.3cm{
          \one\ar[dr]\ar@{.}[dd]&&&&&&\two\ar@{.}[dd]\\
          &\three\ar[dr]&&\four\ar[dr]&&\five\ar[dr]\ar[ur]&\\
          \six\ar[ur]\ar[dr]\ar@{.}[dd]&&\seven\ar[ur]\ar[dr]\ar[r]&\eight\ar[r]&\nine\ar[ur]\ar[dr]&&\ten\ar@{.}[dd]\\
          &\eleven\ar[ur]&&\twive\ar[ur]&&\thirteen\ar[ur]\ar[dr]&\\
          \fourteen\ar[ur]&&&&&&\fifteen
           }\]

\vskip0.2in

Note that there are exactly 14 basis maximal rigid
$\Lambda_{A_3}$-modules up to isomorphism. We list non projective
direct summands of every maximal rigid $\Lambda_{A_3}$-modules as
follows.

$R_1=\six\oplus\three\oplus\eleven$, \qquad$R_2=\seven\oplus\three\oplus\eleven$,

\vskip0.2in
$R_3=\seven\oplus\three\oplus\twive$, \qquad$R_4=\six\oplus\thirteen\oplus\eleven$,

\vskip0.2in
$R_5=\four\oplus\thirteen\oplus\eleven$, \qquad$R_6=\four\oplus\thirteen\oplus\nine$,

\vskip0.2in
$R_7=\four\oplus\twive\oplus\nine$, \qquad$R_8=\four\oplus\twive\oplus\seven$,

\vskip0.2in
$R_9=\five\oplus\thirteen\oplus\nine$, \qquad$R_{10}=\twive\oplus\thirteen\oplus\nine$,

\vskip0.2in
$R_{11}=\twive\oplus\thirteen\oplus\three$, \qquad$R_{12}=\four\oplus\seven\oplus\eleven$,

\vskip0.2in
$R_{13}=\six\oplus\thirteen\oplus\five$, \qquad$R_{14}=\six\oplus\three\oplus\five$.

\vskip0.2in

The mutation graph of basic maximal rigid $\Lambda_{A_3}$-modules is
 following.

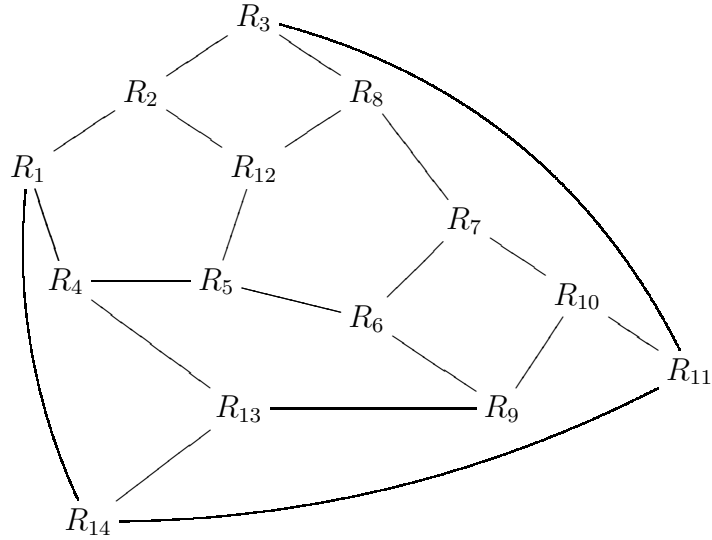
\[\xy
(2,0)*+{R_1}="R_1", +<1.5cm,1cm>*+{R_2}="R_2", +<1.5cm,-1cm>*+{R_{12}}
="R_12",+<0cm,2cm>*+{R_3}="R_3",
+<1.5cm,-1cm>*+{R_8}="R_8",+<-2cm,-2.5cm>*+{R_5}=
"R_5",+<-2cm,0cm>*+{R_4}="R_4",+<+4cm,-0.5cm>*+{R_6}="R_6"
+<1.3cm,1.3cm>*+{R_7}="R_7",+<1.5cm,-1cm>*+{R_{10}}
="R_10",+<-1cm,-1.5cm>*+{R_{9}}="R_9",+<-3.5cm,0cm>*+{R_{13}}="R_13"
+<-2cm,-1.5cm>*+{R_{14}}="R_14",+<8cm,2cm>*+{R_{11}}="R_11"

\POS"R_1" \ar@{-} "R_2"

\POS"R_2" \ar@{-} "R_12"

 \POS"R_2" \ar@{-} "R_3"

 \POS"R_3" \ar@{-} "R_8"

 \POS"R_12" \ar@{-} "R_8"

 \POS"R_12" \ar@{-} "R_5"

 \POS"R_4" \ar@{-} "R_5"

 \POS"R_4" \ar@{-} "R_1"

 \POS"R_6" \ar@{-} "R_5"

 \POS"R_4" \ar@{-} "R_1"

 \POS"R_6" \ar@{-} "R_5"

 \POS"R_6" \ar@{-} "R_7"

 \POS"R_8" \ar@{-} "R_7"

 \POS"R_7" \ar@{-} "R_10"

 \POS"R_9" \ar@{-} "R_10"

 \POS"R_9" \ar@{-} "R_6"

 \POS"R_9" \ar@{-} "R_13"

 \POS"R_13" \ar@{-} "R_4"

 \POS"R_13" \ar@{-} "R_14"

 \POS"R_1" \ar@{-}@/_2ex/ "R_14"

 \POS"R_11" \ar@{-}@/^3ex/ "R_14"

 \POS"R_11" \ar@{-} "R_10"

 \POS"R_11" \ar@{-}@/_5ex/ "R_3"

\endxy \]

\vskip0.2in

According to Corollary 3.10, this picture is also the tilting graphs
for endomorphism algebras of all maximal rigid
$\Lambda_{A_3}$-modules, and every such endomorphism algebras has 14
basic tilting modules up to isomorphism.

\vskip0.2in

{\bf Remarks.}\ We conjecture that Theorem 3.7 is also true for
preprojective algebras of tame representation type. In this case,
the AR-quivers of the preprojective algebras are of tubular type.

\section{The connectedness of mutation graphs of maximal rigid modules}

In this section, we investigate  the connectedness of mutation
graphs of maximal rigid modules over preprojective algebras of
representation finite type or tame type and prove Theorem 3 promised
in the introduction.

\vskip0.2in

It is well known that a preprojective algebra $\Lambda$ is of tame
type if and only if it is of type $A_5$ and $D_4$.  In this case,
their AR-quivers are of tubular type which are the following.

\vskip0.2in

We denote by $\Lambda_5$  the preprojective algebra of type $A_5$,
then the ordinary quiver of $\Lambda_5$ is
\[\xymatrix{ \overline{A}_5: 1\ar@<1ex>[r]^{\alpha_1}&2\ar@<1ex>[r]^{\alpha_2}
\ar@<1ex>[l]^{\alpha_1^*}&3\ar@<1ex>[r]^{\alpha_3}\ar@<1ex>[l]^{\alpha_2^*}&
            4\ar@<1ex>[r]^{\alpha_4}\ar@<1ex>[l]^{\alpha_3^*}&5\ar@<1ex>[l]^{\alpha_4^*}}, \]
and  $\Lambda_5= k\overline{A}_5/I$ with $I$ generated by relations
$\{\alpha_1\alpha_1^*, \alpha_1^*\alpha_1+\alpha_2\alpha_2^*,
\alpha_2^*\alpha_2+\alpha_3\alpha_3^*,
\alpha_3^*\alpha_3+\alpha_4\alpha_4^*, \alpha_4^*\alpha_4\}$.

\vskip0.2in

Note that $\Lambda_5$ admits a Galois covering
$\widetilde{\Lambda_5}$:
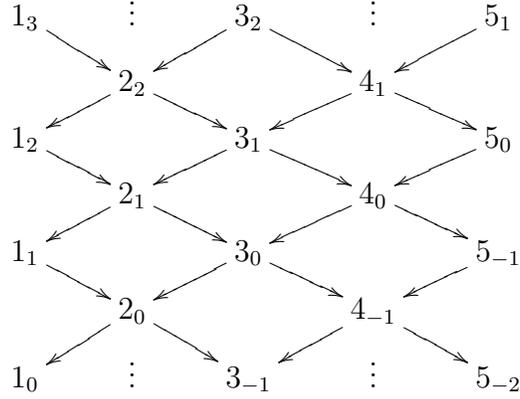
\[\xymatrix@R=0.2cm{1_3\ar[rd]&\vdots&3_2\ar[ld]\ar[rd]&\vdots&5_1\ar[ld] \\
          &2_2\ar[ld]\ar[rd]&&4_1\ar[ld]\ar[rd]& \\
          1_2\ar[rd]&&3_1\ar[rd]\ar[ld]&&5_0\ar[ld] \\
          &2_1\ar[ld]\ar[rd]&&4_0\ar[ld]\ar[rd]& \\
          1_1\ar[rd]&&3_0\ar[ld]\ar[rd]&&5_{-1}\ar[ld] \\
          &2_0\ar[ld]\ar[rd]&&4_{-1}\ar[ld]\ar[rd]& \\
          1_0&\vdots&3_{-1}&\vdots&5_{-2}\\
          }\]
with the mesh relations and zero relations. All $\Lambda_5$-module
can be obtained by applying the push down functor to the
$\widetilde{\Lambda_5}$-modules, and  $\widetilde{\Lambda_5}$ can be
regarded as the repetitive algebras of the tubular algebra $\Delta$:
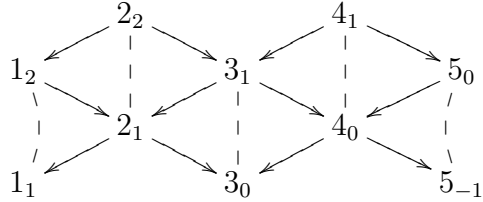
\[\xymatrix@R=0.2cm{
          &2_2\ar[ld]\ar[rd]\ar@{--}[dd]&&4_1\ar[ld]\ar[rd]\ar@{--}[dd]& \\
          1_2\ar[rd]\ar@{--}@/^/[dd]&&3_1\ar[rd]\ar[ld]\ar@{--}[dd]&&5_0\ar[ld]\ar@{--}@/_/[dd] \\
          &2_1\ar[ld]\ar[rd]&&4_0\ar[ld]\ar[rd]& \\
          1_1&&3_0&&5_{-1}
          }\]
of tubular type (6,3,2). We have $\underline{\Mod}\
\widetilde{\Lambda_5}\cong \mathcal{D}^b(\Mod\Delta)\cong
\mathcal{D}^b(coh(\mathbb{X}))$ by the theorems of Happel, Geigle
and Lenzing, where $\mathbb{X}$ is a weighted projective line of
type (6,3,2),  see \cite[setion9 and section 19]{GLS2} for details.

\vskip0.2in

Let $\Lambda_{D_4}$ be the preprojective algebra of type $D_4$, then
the  ordinary quiver of $\Lambda_{D_4}$ is
\[\xymatrix{\overline{D}_4: 2\ar@<1ex>[r]^{\alpha_1^*}&1\ar@<1ex>[l]^{\alpha_1}
\ar@<1ex>[d]^{\alpha_2}\ar@<1ex>[r]^{\alpha_3}&
4\ar@<1ex>[l]^{\alpha_3^*}\\&3\ar@<1ex>[u]^{\alpha_2^*}&}, \] and
$\Lambda_{D_4}= k\overline{D}_4/I$ with $I$ generated by relations
$\{\alpha_1^*\alpha_1, \alpha_2^*\alpha_2, \alpha_3^*\alpha_3,
\alpha_1\alpha_1^*+\alpha_2\alpha_2^*+\alpha_3\alpha_3^*\}$.

\vskip0.2in

Note that  $\Lambda_{D_4}$  has a Galois covering
$\widetilde{\Lambda_{D_4}}$ as follows:
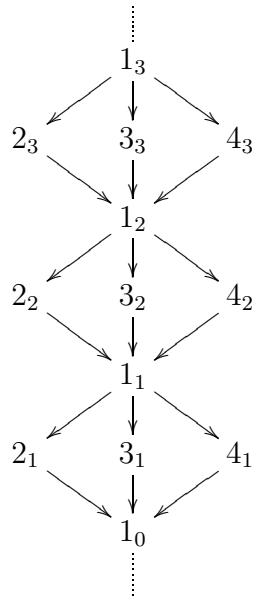
\[\xymatrix@R=0.5cm{  &\ar@{..}[d]&\\
           &1_3\ar[ld]\ar[d]\ar[rd]&\\
           2_3\ar[rd]&3_3\ar[d]&4_3\ar[ld]\\
           &1_2\ar[ld]\ar[d]\ar[rd]&\\
           2_2\ar[rd]&3_2\ar[d]&4_2\ar[ld]\\
           &1_1\ar[ld]\ar[d]\ar[rd]&\\
           2_1\ar[rd]&3_1\ar[d]&4_1\ar[ld]\\
           &1_0\ar@{..}[d]&\\
            &&}\]
with the mesh relations and the zero relations. It can be regarded
as the repetitive algebra $\Delta$:
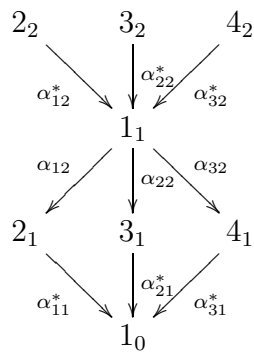
\[\xymatrix{2_2\ar[rd]_{\alpha_{12}^*}&3_2\ar[d]^{\alpha_{22}^*}&4_2\ar[ld]^{\alpha_{32}^*}\\
           &1_1\ar[ld]_{\alpha_{12}}\ar[d]^{\alpha_{22}}\ar[rd]^{\alpha_{32}}&\\
           2_1\ar[rd]_{\alpha_{11}^*}&3_1\ar[d]^{\alpha_{21}^*}&4_1\ar[ld]^{\alpha_{31}^*}\\
           &1_0&}\]
the relations are $\{\alpha_{12}^*\alpha_{12},
\alpha_{22}^*\alpha_{22}, \alpha_{32}^*\alpha_{32},
\alpha_{12}\alpha_{11}^*+\alpha_{22}\alpha_{21}^*+\alpha_{32}\alpha_{31}^*\}$.

\vskip0.2in

It is a tubular algebra obtained through one point extensions of the
$\widetilde{D_4}$ tame concealed algebra $\Delta^0$:
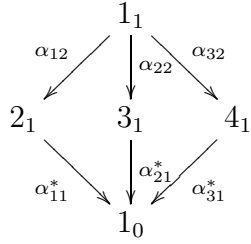
\[\xymatrix{&1_1\ar[ld]_{\alpha_{12}}\ar[d]^{\alpha_{22}}\ar[rd]^{\alpha_{32}}&\\
           2_1\ar[rd]_{\alpha_{11}^*}&3_1\ar[d]^{\alpha_{21}^*}&4_1\ar[ld]^{\alpha_{31}^*}\\
           &1_0&}\]
Thus $\Delta$ is of tubular type (3,3,3). Again,
$\underline{\Mod}\widetilde{\Lambda_{D_4}}\cong
\mathcal{D}^b(\Mod\Delta)\cong\mathcal{D}^b(coh(\mathbb{X}))$, where
$\mathbb{X}$ is a weighted projective line of type (3,3,3), see
[12, section 9 and section 19] for details.

\vskip0.2in

Let $G=\mathbb{Z}$ be the Galois group of the Galois covering
$F:\widetilde{\Lambda}\rightarrow \Lambda$. $G$ has an action on the
$\widetilde{\Lambda}$-modules $X$, that is for every vector space
$X_{k_j}$ of $X$ corresponding to the vertex $k_j$, we get $X^{(i)}$
with $X_{k_{j+i}}=X_{k_j}$ and keep the maps between the vector
spaces. Let $F$ be the push down functor from $\Mod\
\widetilde{\Lambda}$ to $\Mod\ \Lambda$. Then we have
$\Hom_{\Lambda}(F(X),F(X))=\sum_{i\in\mathbb{Z}}\Hom_{\widetilde{\Lambda}}(X,X^{(i)})$.

\vskip0.2in

Let $\mathcal{C}$ be the cluster category of a hereditary abelian
category with cluster-tilted objects in the sense of \cite{BMRRT1,
ZB}. According to \cite[Proposition 3.5]{BMRRT1},   the tilting
graph of cluster-tilted objects in $\mathcal{C}$ is connected if
$\mathcal{C}$ is the cluster category of a finite dimensional
hereditary algebra.

\vskip0.2in

{\bf Theorem 4.1.}\  {\it Let $\Lambda$ be a preprojective algebra
of finite or tame representation type. Then the mutation graph of
basic maximal rigid $\Lambda$-modules is connected.}

\vskip0.1in

{\bf Proof.} It is well known that $\underline{\Mod}\ \Lambda$ is
2-Calabi-Yau. And it is clear that the basic maximal rigid modules
of $\Mod\ \Lambda$ are in bijection with the basic cluster-tilted
objects in $\underline{\Mod}\ \Lambda$ and the mutation graphs of
them are the same by definitions. So we only need to consider the
mutation graph of the basic cluster-tilted objects in
$\underline{\Mod}\ \Lambda$.

If $\Lambda$ is of type $A_2$, $A_3$ or $A_4$, then the AR-quiver of
$\underline{\Mod}\ \Lambda$ is the same with the quivers of the
cluster category $\mathcal {C}$ of $A_1$, $A_3$, and $D_6$
respectively, see \cite[section 20.1]{GLS2}. Hence, the
cluster-tilted objects in $\underline{\Mod}\ \Lambda$ are in
bijection with the cluster-tilted objects in $\mathcal{C}$. Hence
the mutation graph of the basic cluster-tilted objects in
$\underline{\Mod}\ \Lambda$ is connected by \cite[Proposition
3.5]{BMRRT1}.

If $\Lambda$ is of type $A_5$, we know that
$\underline{\Mod}\ \widetilde{\Lambda_5}\cong
\mathcal{D}^b(\Mod\Delta)\cong \mathcal{D}^b(coh(\mathbb{X}))$, and
by \cite[section 14.5,14.6]{GLS2} we know that $\underline{\Mod}\
{\Lambda_5}$ is a fundamental domain of $\underline{\Mod}\
\widetilde{\Lambda_5}$ under the action of the Galois group
$\mathbb{Z}$. So, $\underline{\Mod}\ {\Lambda_5}\cong
\underline{\Mod}\ \widetilde{\Lambda_5}/(1)$ as the orbit category,
where $(1)$ is the generator of the Galois group. The cluster
category of $\mathcal{D}^b(coh(\mathbb{X}))$ is by definition
$\mathcal
{C}=\mathcal{D}^b(coh(\mathbb{X}))/\tau^{-1}[1]\cong\underline{\Mod}\
\widetilde{\Lambda_5}/\tau^{-1}[1]$, where $\tau^{-1}$ is the inverse
of the AR translation and [1] is the shift functor. By \cite[Lemma
6.1]{GS}, $(1)\cong\tau^{-1}[1]$, so we have $\mathcal
{C}\cong\underline{\Mod}\ {\Lambda_5}$. By \cite[Theorem 8.8]{BKL},
the tilting graph of $\mathcal{C}$ is connected, so the mutation
graph of $\underline{\Mod}\ {\Lambda_5}$ is connected. The $D_4$
case can be proved similarly.    This completes the proof.
$\hfill\Box$

\vskip0.5in

{\bf Acknowledgement.} The authors would like to thank Professor
C.M.Ringel for many useful comments and helpful discussions.

\end{document}